\documentclass[a4paper, 12pt] {article}
\usepackage[cp1251]{inputenc}
\usepackage[english,russian]{babel}
\usepackage{mathtext}
\usepackage{amsfonts}
\usepackage{cite}

\begin{document}

\noindent УДК 517.9

\vspace*{5mm}

\begin{center}
{\Large\bf Об успокоении системы управления произвольного порядка с глобальным последействием на дереве} \footnote{Работа выполнена при
поддержке РНФ, проект №~22-21-00509, https://rscf.ru/project/22-21-00509/}
\end{center}

\begin{center}
{\large {\bf С.\,А.~Бутерин}\\[3mm]
Саратовский университет, г.~Саратов\\[3mm]
 buterinsa@sgu.ru}
\end{center}

{\bf Аннотация.} Исследуется задача об успокоении управляемой системы, описываемой функ\-ционально-дифференциальными уравнениями натурального
порядка $n$ нейт\-рального типа с негладкими комплексными коэффициентами на произвольном дереве с глобальным запаздыванием. Последнее
означает, что запаздывание распространяется через внутренние вершины дерева. Минимизация функционала энергии системы приводит к вариационной
задаче. Установлена ее эквивалентность некоторой самосопряженной краевой задаче на дереве для уравнений порядка $2n$ с нелокальными
квазипроизводными и разнонаправленными сдвигами аргумента, а также условиями типа Кирхгофа, возникающими во внутренних вершинах. Доказана
однозначная разрешимость обеих задач.

\medskip
{\it Ключевые слова}: квантовый граф, функционально-дифференциальное уравнение, глобальное запаздывание, задача оптимального управления,
вариационная задача, нелокальная квазипроизводная, временной граф



\begin{center}
{\bf 1. Введение}
\end{center}

Дифференциальные операторы на геометрических графах, часто называемые квантовыми графами, активно исследуются с прошлого века в связи с
многочисленными приложениями \cite{Mont-70, PPPB, BerkKuch}. Такие операторы возникают при исследовании процессов в сложных системах,
представимых в виде пространственных сетей, то есть наборов одномерных континуумов, взаимодействующих только через концы\cite{PPPB}. Примером
являются упругие струнные сетки, в узлах которых помимо условий непрерывности характерными являются условия Кирхгофа, выражающие баланс
натяжений.

В настоящей работе переменная на ребрах графа отождествляется со временем, когда в каждой внутренней вершине процесс разветвляется на
несколько параллельных процессов по числу выходящих из нее ребер. В этом случае тоже могут возникать условия Кирхгофа. А именно, им будет
удовлетворять такая траектория течения процесса, которая является оптимальной с учетом сразу всех перспектив.

На временн\' ых графах становится актуальным рассмотрение также процессов с последействием. Однако до недавнего времени работы, посвященные
нелокальным операторам на графах, относились в основном к {\it локально} нелокальному случаю, когда соответствующее нелокальное уравнение на
каждом ребре решается независимо от уравнений на остальных ребрах \cite{Nizh-12, Bon18-1, HuBondShYan19, Hu20, WangYang-21}. При этом не
освещался вопрос, как можно было бы определить оператор в {\it глобально} нелокальном случае и, в частности, как описать процесс с глобальным
последействием на всем графе.

Чтобы восполнить этот пробел, в \cite{But23-RM} была предложена концепция функционально-дифферен\-циальных операторов на графах с глобальным
запаздыванием, которое ''проходит'' через вершины графа, а также исследовались обратные спектральные задачи для таких операторов (см. тж.
\cite{But23-M}). С помощью этой концепции в \cite{But23-arXiv} на графы был распространен и другой класс задач, для которого уже естественно
отождествление графа со временем. Речь идет об успокоении управляемой системы с последействием для уравнения первого порядка запаздывающего
типа.

На интервале эта задача впервые была поставлена и исследована Н.Н.~Красовским~\cite{Kras-68} для уравнений с постоянными вещественными
коэффициентами. Позднее А.Л. Скубачевский \cite{Skub-94} рассмотрел ее обобщение на случай, когда уравнение содержит также старшие члены с
запаздыванием, то есть имеет нейтральный тип. Минимизация функционала энергии приводит к соответствующей вариационной задаче. Была доказана
ее эквивалентность некоторой самосопряженной краевой задаче для уравнения второго порядка с двунаправленными сдвигами аргумента, и
установлена однозначная разрешимость этой краевой задачи. При этом нейтральный тип исходного уравнения приводит к понятию обобщенного решения
соответствующей краевой задачи, чья первая производная может терять гладкость. В \cite{Skub-94} получены также необходимые и достаточные
условия гладкости, а в монографии \cite{Skub-97} приводится общая теория таких решений. Указанные результаты обобщены в \cite{AdSkub-19,
AdSkub-20} на случай управляемой системы, заданной уравнениями нейтрального типа с непрерывно дифференцируемыми коэффициентами в пространстве
вектор-функций. Отдельно следует отметить аналогичные результаты для случая, когда запаздывание не постоянно, а является пропорциональным
времени сжатием \cite{Ross-14}.

В \cite{But23-arXiv} данная задача при запаздывающем типе исходных уравнений была распространена на произвольное дерево. Полученную систему
управления можно интерпретировать следующим образом. В каждый момент времени, соответствующий какой-либо внутренней вершине дерева,
появляется несколько вариантов дальнейшего течения процесса (по числу выходящих из этой вершины ребер). Требуется подобрать управление,
приводящее систему в состояние равновесия, независимо от того, какой именно набор сценариев в итоге будет реализован, то есть с учетом сразу
всех перспектив. Были установлены существование и единственность оптимальной траектории, которая, как и выяснилось, дополнительно
удовлетворяет условиям типа Кирхгофа во всех внутренних вершинах дерева (см. ниже замечание~1).

В настоящей работе мы переходим к уравнениям нейтрального типа на графах. Однако рассматривается более общий случай управляемой системы,
уравнения которой имеют произвольный натуральный порядок~$n,$ а коэффициенты являются комплекснозначными функциями из $L_2,$ тогда как самые
старшие принадлежат $L_\infty.$

Для произвольного дерева установлено, что соответствующая вариационная задача эквивалентна некоторой самосопряженной краевой задаче для
уравнений порядка $2n$ с разнонаправленными сдвигами аргумента и дополнительными $n$ условиями типа Кирхгофа в каждой внутренней вершине
(теорема~5). При этом дальнейшее развитие получает концепция обобщенного решения краевой задачи, которое перестает быть прерогативой
исключительно нейтрального типа исходных уравнений даже в случае постоянных коэффициентов. Теперь для ее постановки требуется введение
семейства нелокальных квазипроизводных порядков от $n$ до $2n,$ параллельно позволяющих охватить и рассматриваемый случай негладких
коэффициентов. При этих же условиях доказана однозначная разрешимость данной краевой задачи.

\medskip
Поскольку полученные результаты являются новыми даже в случае интервала, в следующем разделе они формулируются отдельно для этого
случая. Там же демонстрируется преемственность понимания решения краевой задачи в терминах введенных квазипроизводных с понятием обобщенного
решения в \cite{Skub-94}. Для иллюстрации специфики графа в разделе~3 результаты приводятся отдельно для графа-звезды. Постановка
вариационной задачи на произвольном дереве дается в разделе~4. Раздел~5 посвящен доказательству ее эквивалентности соответствующей краевой
задаче, чья однозначная разрешимость установлена в разделе~6. В последнем разделе проводится сравнение введенных квазипроизводных в локальном
случае с классическими квазипроизводными для обыкновенных дифференциальных операторов.

\newpage

\begin{center}
{\bf 2. Случай интервала}
\end{center}

Рассмотрим управляемую систему, описываемую уравнением нейтрального типа
\begin{equation}\label{1.1}
\ell y(t):=\sum_{k=0}^n\Big(b_k(t) y^{(k)}(t)+c_k(t) y^{(k)}(t-\tau)\Big)=u(t), \quad t>0,
\end{equation}
с постоянным запаздыванием $\tau>0,$ где коэффициенты $b_k=b_k(t),$ $c_k=c_k(t),$ а также управление $u=u(t)$ предполагаются
комплекснозначными функциями, причем
\begin{equation}\label{1.1-1}
\forall\, a>0 \quad b_n,\,\frac1{b_n},\,c_n\in L_\infty(0,a), \quad u,\,b_k,\,c_k\in L_2(0,a), \;\; k=\overline{0,n-1}.
\end{equation}
Предыстория системы описывается условием
\begin{equation}\label{1.2}
y(t)=\varphi(t), \quad t\in(-\tau,0),
\end{equation}
с заданной комплекснозначной функцией $\varphi(t)\in W_2^n[-\tau,0],$ а также условиями
\begin{equation}\label{1.2-1}
y^{(k)}(0)=\varphi^{(k)}(0), \quad k=\overline{0,n-1}.
\end{equation}

Приводя уравнение (\ref{1.1}) к системе первого порядка и последовательно применяя теорему~1 из гл.~V в \cite{Neum} для каждого интервала
$((j-1)\tau,j\tau),\;j\in{\mathbb N},$ легко показать, что при всех $a>0$ задача Коши (\ref{1.1})--(\ref{1.2-1}) имеет единственное решение
$y\in W_2^n[0,a].$

Зафиксируем $T>2\tau$ и рассмотрим задачу об успокоении системы (\ref{1.1})--(\ref{1.2-1}), которое характеризуется условием $y(t)=0$ при
$t\ge T.$ Для этого достаточно применить управление $u(t)\in L_2(0,T),$ приводящее систему в состояние равновесия
\begin{equation}\label{1.3}
y(t)=0, \quad t\in[T-\tau,T],
\end{equation}
а затем убрать воздействие на нее, положив $u(t)=0$ при $t>T.$

Поскольку такое $u(t)$ на $(0,T)$ не единственно, естественной является попытка минимизировать затрачиваемые при этом усилия
$\|u\|_{L_2(0,T)},$ которая приводит к вариационной задаче о минимуме соответствующего функционала энергии
\begin{equation}\label{1.4}
J(y)=\int_0^T |\ell y(t)|^2\,dt \to\min
\end{equation}
при условиях (\ref{1.2})--(\ref{1.3}). Обозначим
\begin{equation}\label{1.4-00}
\tilde \ell_{k}y(t):= \overline{b_k(t)}\ell y(t) +\overline{c_k(t+\tau)}\ell y(t+\tau), \quad 0< t < T-\tau, \quad k=\overline{0,n},
\end{equation}
и введем квазипроизводные
\begin{equation}\label{1.4-0}
y^{\langle n\rangle}(t):=\tilde\ell_{n} y(t), \quad y^{\langle n+l\rangle}(t):=\tilde\ell_{n-l} y(t) -(y^{\langle n+l-1\rangle})'(t), \quad
l=\overline{1,n}.
\end{equation}

\medskip
{\bf Теорема 1. }{\it Функция $y(t)\in W_2^n[-\tau,T]$ является решением вариационной задачи (\ref{1.2})--(\ref{1.4}) тогда и только тогда,
когда она удовлетворяет условиям
\begin{equation}\label{1.4-1}
y^{\langle k\rangle}(t)\in W_1^1[0,T-\tau], \quad k=\overline{n,2n-1},
\end{equation}
и является решением самосопряженной краевой задачи ${\mathcal B}$ для уравнения
$$
y^{\langle 2n\rangle}(t)=0, \quad 0<t<T-\tau,
$$
при условиях (\ref{1.2})--(\ref{1.3}).}

\newpage
Вместе с теоремой~1 следующая теорема дает существование и единственность решения вариационной задачи (\ref{1.2})--(\ref{1.4}).

\medskip
{\bf Теорема 2. }{\it Краевая задача ${\cal B}$ имеет единственное решение $y(t)\in W_2^n[-\tau,T],$ удовлетворяющее условиям (\ref{1.4-1}).
При этом справедлива оценка
$$
\|y\|_{W_2^n[-\tau,T]} \le C\|\varphi\|_{W_2^n[-\tau,0]},
$$
где $C$ не зависит от $\varphi(t).$}

\medskip
Теоремы~1,~2 можно получить как следствия теорем~3,~4 для графа-звезды, в свою очередь являющихся частными случаями теорем~5,~6 для
произвольного дерева.

\medskip
Укажем на преемственность понимания решения задачи ${\cal B}$ в терминах введенных квазипроизводных с обобщенным решением в смысле
\cite{Skub-94}. Для этого положим
\begin{equation}\label{1.4-3}
n=1,\quad b_1=1,\quad c_1=a,\quad b_0=b,\quad c_0=c,
\end{equation}
где $a,b,c\in {\mathbb R}.$ Тогда уравнение (\ref{1.1}) становится уравнением из \cite{Skub-94}:
\begin{equation}\label{1.4-4}
\ell y(t)=y'(t)+ay'(t-\tau)+ by(t) +cy(t-\tau)=u(t).
\end{equation}
При этом выражения в (\ref{1.4-00}) примут вид
$$
\tilde\ell_0 y(t)=b\ell y(t) +c\ell y(t+\tau), \quad \tilde\ell_1 y(t)=\ell y(t) +a\ell y(t+\tau).
$$
Согласно (\ref{1.4-0}) и (\ref{1.4-4}), будем иметь
$$
y^{\langle 1\rangle}(t)=Ay(t)+By(t), \quad y^{\langle 2\rangle}(t)=\tilde\ell_0 y(t) -(Ay)'(t)-(By)'(t),
$$
где
\begin{equation}\label{1.4-5}
Ay(t)= (1+a^2)y'(t) +ay'(t-\tau)+ay'(t+\tau),
\end{equation}
$$
By(t)=(ac+b)y(t)+cy(t-\tau)+aby(t+\tau).
$$
Далее, вычисляем $y^{\langle 2\rangle}(t)= -(Ay)'(t) +Cy(t),$ где
$$
Cy(t)= (ab-c)(y'(t-\tau) -y'(t+\tau)) +(b^2+c^2)y(t) +bc(y(t-\tau) +y(t+\tau)).
$$
Таким образом, краевая задача ${\cal B}$ в частном случае (\ref{1.4-3}) совпадает с краевой задачей (2), (3), (8) в \cite{Skub-94}.

Поскольку $By(t)$ всегда принадлежит $ W_2^1[0,T-\tau],$  условие (\ref{1.4-1}) равносильно $Ay(t)\in W_1^1[0,T-\tau].$ Последнее, в силу
$Cy(t)\in L_2(0,T-\tau),$ равносильно, в свою очередь, принадлежности $y$ классу функций, в смысле которого обобщенные решения понимаются в
\cite{Skub-94}, а именно: $Ay(t)\in W_2^1[0,T-\tau],$ коль скоро $y^{\langle 2\rangle}(t)=0.$

Итак, решение задачи ${\cal B}$ в частном случае (\ref{1.4-3}), понимаемое в смысле абсолютной непрерывности $y^{\langle 1\rangle},$
соответствует концепции обобщенного решения из \cite{Skub-94}.

\medskip
Как отмечено в \cite{Skub-94}, а также непосредственно вытекает из представления (\ref{1.4-5}), рассмотрение обобщенных решений в случае
$n=1$ для постоянных коэффициентов актуально только при $a\ne0,$ то есть при нейтральном типе уравнения (\ref{1.4-4}).

Покажем, что при $n>1$ обобщенное решение в том смысле, что какая-либо его производная порядка с $n$ по $2n-1$ теряет гладкость, может
возникать и при запаздывающем типе уравнения (\ref{1.1}) даже в случае постоянных коэффициентов.

Для этой цели рассмотрим случай, когда
\begin{equation}\label{1.4-6}
n=2, \quad b_2=c_1=1, \quad b_0=b_1=c_0=c_2=0,
\end{equation}
то есть
$$
\ell y(t)=y''(t)+y'(t-\tau).
$$
В силу (\ref{1.4-00}) имеем
$$
\tilde\ell_0y(t)=0, \quad \tilde\ell_1y(t)=\ell y(t+\tau), \quad \tilde\ell_2y(t)=\ell y(t),
$$
что вместе с  (\ref{1.4-0}) приводит к квазипроизводным
$$
y^{\langle 2\rangle}(t)=\tilde\ell_2 y(t)=y''(t)+y'(t-\tau),
$$
$$
y^{\langle 3\rangle}(t)=\tilde\ell_1 y(t)-(y^{\langle 2\rangle})'(t) =-y'''(t)+y''(t+\tau)-y''(t-\tau)+y'(t).
$$
Если, например, $T=4\tau$ и $\varphi(t)\in W_2^2[-\tau,0]\setminus W_1^3[-\tau,0],$ то (\ref{1.2}) и (\ref{1.2-1}) дают
$$
y'''(t)\in W_1^1[0,3\tau] \;\;\Rightarrow\;\; y^{\langle 3\rangle}(t)\notin W_1^1[0,3\tau].
$$
Заметим, что последняя импликация будет справедливой и при $\varphi(t)\in W_1^3[-\tau,0],$ если потребовать $\varphi''(0)\ne y''(0),$ что не
исключается условиями (\ref{1.2-1}) при $n=2.$

Таким образом, мы показали, что обычное, то есть ''не обобщенное'', решение $y\in W_2^2[-\tau,T]\cap W_1^4[0,T-\tau]$ краевой задачи ${\cal
B}$ в случае (\ref{1.4-6}) может не существовать, поскольку, в силу теоремы~2, ее решение всегда удовлетворяет условиям (\ref{1.4-1}).

\begin{center}
{\bf 3. Граф типа звезды}
\end{center}

Предположим, что до момента времени $t=T_1,$ ассоциированного с единственной внутренней вершиной $v_1$ графа $\Gamma_m,$ изображенного на
рис.~1, наша управляемая система с запаздыванием $\tau<T_1$ на $\Gamma_m$ описывается уравнением
\begin{equation}\label{1.5}
\ell_1 y(t):=\sum_{k=0}^n\Big(b_{k,1}(t)y_1^{(k)}(t)+c_{k,1}(t)y_1^{(k)}(t-\tau)\Big)=u_1(t), \quad 0<t<T_1,
\end{equation}
заданным на ребре $e_1,$ и имеет предысторию, определяемую соотношением
\begin{equation}\label{1.9}
y_1(t)=\varphi(t)\in W_2^n[-\tau,0], \quad t\in(-\tau,0),
\end{equation}
и начальными условиями
\begin{equation}\label{1.9-1}
y_1^{(k)}(0)=\varphi^{(k)}(0), \quad k=\overline{0,n-1}.
\end{equation}

\begin{center}
\unitlength=0.9mm
\begin{picture}(80,100)
\multiput(-10,52)(-1,0){14}{\circle*{0.5}} \put(-23.4,52){\circle*{1}} \put(-29,54){\small $-\tau$}

 \put(-10,52){\line(1,0){50}}    \put(11,58){$e_1$} \put(-10,52){\circle*{1}}
 \put(-11,47){\small $v_0$}    \put(-11,54){\small $0$}
                                                              \put(39,58){\small $0$}
                                                              \put(43.5,56){\small $0$}
                                                              \put(43,45){\small $0$}
 \put(26.6,52){\circle*{1}}                                   \put(14,46){\small $T_1-\tau$}
 \put(40,52){\circle*{1}}        \put(35,47){\small $v_1$}     \put(34,54){\small $T_1$}

 \put(40,52){\vector(1,3){15}}     \put(55,98){\small $e_2$} \put(47,73){\circle*{1}} \put(48,68.5){\small $T_2-\tau$}
 \put(51.3,85.9){\circle*{1}}   \put(45,86){\small $T_2$}

 \put(40,52){\vector(2,1){61}}     \put(102.3,82.3){\small $e_3$} \put(78,71){\circle*{1}}
 \put(78.5,66.5){\small $T_3-\tau$} \put(90,77){\circle*{1}} \put(86,80){\small $T_3$}

 \multiput(40,52)(2,-1){22}{\circle*{0.5}}

 \put(40,52){\vector(1,-3){13}}     \put(52.5,9){\small $e_m$} \put(45,37){\circle*{1}}
 \put(29,33){\small $T_m-\tau$} \put(49.3,24.1){\circle*{1}} \put(51.5,24){\small $T_m$}

\put(22,0){\small Рис. 1. Граф $\Gamma_m$}
\end{picture}
\end{center}

Далее, при $t=T_1,$ то есть в вершине $v_1,$ процесс разветвляется на $m-1$ (по числу оставшихся ребер) независимых параллельных процессов,
заданных уравнениями
\begin{equation}\label{1.6}
\ell_j y(t):=\sum_{k=0}^n\Big(b_{k,j}(t)y_j^{(k)}(t)+c_{k,j}(t)y_j^{(k)}(t-\tau)\Big)=u_j(t), \quad t>0, \quad j=\overline{2,m},
\end{equation}
но имеющих общую историю, определяемую уравнением (\ref{1.5}) с предысторией (\ref{1.9}) и (\ref{1.9-1}), а также условиями прохождения
запаздывания через вершину $v_1:$
\begin{equation}\label{1.7}
y_j(t)=y_1(t+T_1), \quad t\in(-\tau,0), \quad j=\overline{2,m}.
\end{equation}
Соответственно $j\!$-ое уравнение в (\ref{1.6}) задано на ребре $e_j$ графа $\Gamma_m,$ представляющем собой, вообще говоря, бесконечный луч,
выходящий из $v_1.$

Кроме того, естественно наложить условия непрерывности в вершине $v_1:$
\begin{equation}\label{1.8}
y_j^{(k)}(0)=y_1^{(k)}(T_1), \quad k=\overline{0,n-1}, \quad  j=\overline{2,m}.
\end{equation}

Зафиксируем произвольные числа $T_j>\tau,$ $j=\overline{2,m}.$ Успокоение получившейся системы (\ref{1.5})--(\ref{1.8}) будет означать, что
решение $y_j(t)$ уравнения в (\ref{1.6}) для каждого $j=\overline{2,m}$ становится тождественным нулем при $t\ge T_j.$

Другими словами, система должна прийти в положение равновесия заведомо при любом возможном сценарии, допускаемом с момента времени $t=T_1.$

Для этой цели достаточно отыскать управления $u_j(t)\in L_2(0,T_j),$ $j=\overline{1,m},$ приводящие систему в состояние
\begin{equation}\label{1.10}
y_j(t)=0, \quad t\in[T_j-\tau,T_j], \quad j=\overline{2,m}.
\end{equation}

Поскольку набор таких $u_j(t)$ не единственен, разумно искать его из условия минимума функционала энергии. Таким образом, приходим к
вариационной задаче
$$
\sum_{j=1}^m\int_0^{T_j}|\ell_j y_j(t)|^2\,dt\to\min
$$
при условиях (\ref{1.9}), (\ref{1.9-1}) и (\ref{1.7})--(\ref{1.10}), которую для краткости обозначим ${\cal V}.$

Предполагается, что все функции, входящие в уравнения (\ref{1.5}) и (\ref{1.6}), а также $\varphi(t)$ в (\ref{1.9}) комплекснозначны, причем
\begin{equation}\label{1.8-1}
b_{n,j},\frac1{b_{n,j}},c_{n,j}\in L_\infty(0,T_j), \quad u_j, b_{k,j},c_{k,j}\in L_2(0,T_j), \;\; k=\overline{0,n-1}, \quad
j=\overline{1,m}.
\end{equation}

Для $k=\overline{0,n}$ положим
$$
\ell_{k,j} y(t):=\overline{b_{k,j}(t)}\ell_j y(t) + \overline{c_{k,j}(t+\tau)}\ell_j y(t+\tau), \quad 0<t<T_j-\tau, \quad j=\overline{1,m},
$$
$$
\quad\, \ell_{k,1} y(t):=\overline{b_{k,1}(t)}\ell_1 y(t) + \sum_{\nu=2}^m \overline{c_{k,\nu}(t+\tau-T_1)}\ell_\nu y(t+\tau-T_1), \quad
T_1-\tau< t<T_1,
$$
и введем квазипроизводные
\begin{equation}\label{1.8-2}
y_j^{\langle n\rangle}(t):=\ell_{n,j} y(t), \;\; y_j^{\langle n+l\rangle}(t):=\ell_{n-l,j} y(t) -(y_j^{\langle n+l-1\rangle})'(t), \;\;
l=\overline{1,n}, \;\; j=\overline{1,m}.
\end{equation}
Поскольку $y_j^{\langle k\rangle}$ может зависеть не только от $y_j,$ но и от $y_\nu$ на других ребрах~$\Gamma_m,$ то есть при $\nu\ne j,$
квазипроизводные имеют глобально нелокальный характер.

\medskip
{\bf Теорема 3. }{\it Функции $y_1(t)\in W_2^n[-\tau,T_1]$ и $y_j(t)\in W_2^n[0,T_j],$ $j=\overline{2,m},$ образуют решение вариационной
задачи ${\cal V}$ тогда и только тогда, когда они удовлетворяют условиям абсолютной непрерывности соответствующих квазипроизводных:
\begin{equation}\label{1.11-1}
y_1^{\langle k\rangle}(t)\in W_1^1[0,T_1], \quad y_j^{\langle k\rangle}(t)\in W_1^1[0,T_j-\tau], \quad j=\overline{2,m}, \quad
k=\overline{n,2n-1},
\end{equation}
и являются решениями самосопряженной краевой задачи ${\mathcal B}$ для уравнений
$$
y_1^{\langle 2n\rangle}(t)=0, \quad 0<t<T_1, \quad y_j^{\langle 2n\rangle}(t)=0, \quad 0<t<T_j-\tau, \quad j=\overline{2,m},
$$
при условиях (\ref{1.9}), (\ref{1.9-1}), (\ref{1.7})--(\ref{1.10}), а также при дополнительных условиях в $v_1:$
\begin{equation}\label{1.12}
y_1^{\langle k\rangle}(T_1)=\sum_{\nu=2}^m y_\nu^{\langle k\rangle}(0), \quad k=\overline{n,2n-1}.
\end{equation}}

\medskip
{\bf Замечание 1.} Условие (\ref{1.12}) можно охарактеризовать как нелокальное условие типа Кирхгофа. Например, в случае $n=1,$ когда
$b_{1,j}(t)=1$ и $c_{1,j}(t)=0,$ а также $b_{0,j}(t),c_{0,j}(t)\in W_1^1[0,T_j]$ при $j=\overline{1,m},$ абсолютная непрерывность
квазипроизводных $y_j^{\langle 1\rangle}(t)$ равносильна абсолютной непрерывности обычных производных $y_j'(t),$ а условие (\ref{1.12})
принимает вид
$$
y_1'(T_1)+\Big(b_{0,1}(T_1) -\sum_{j=2}^mb_{0,j}(0)\Big) y_1(T_1) +\Big(c_{0,1}(T_1) -\sum_{j=2}^mc_{0,j}(0)\Big) y_1(T_1-\tau)=\sum_{j=2}^m
y_j'(0).
$$
В частности, если выражения в скобках равны нулю, последнее условие является классическим условием Кирхгофа, которое часто возникает как в
общей теории квантовых графов, так и в конкретных приложениях, но относительно пространственной переменной. Например, оно выражает баланс
натяжений в системе~$m$ связанных струн или закон Кирхгофа в электрических цепях. Если же, как в нашем случае, отвечающая графу переменная
ассоциирована со временем, естественно накладывать во внутренних вершинах только условия непрерывности~(\ref{1.8}). Однако, как видно из
теоремы~3, условия типа Кирхгофа возникают и здесь -- а именно, для оптимальной траектории течения процесса с учетом всех допустимых
сценариев.

\medskip
{\bf Теорема 4. }{\it Краевая задача ${\cal B}$ имеет единственное решение
$$
y(t)\in W_2^n[-\tau,T_1], \quad y_j(t)\in W_2^n[0,T_j], \quad j=\overline{2,m},
$$
удовлетворяющее условиям (\ref{1.11-1}). Также найдется $C,$ не зависящее от $\varphi(t),$ при котором справедлива оценка
$$
\|y_1\|_{W_2^n[-\tau,T_1]}+\sum_{j=2}^m \|y_j\|_{W_2^n[0,T_j-\tau]} \le C\|\varphi\|_{W_2^n[-\tau,0]}.
$$}

\medskip
Ниже мы обобщим теоремы~3 и~4 на случай произвольного дерева (теоремы~5 и~6).  Как и $\Gamma_m,$ оно будет содержать изначально бесконечные
ребра, но фактически речь пойдет о компактном дереве ${\cal T},$ получаемом обрезанием бесконечных граничных ребер в соответствующих точках
$T_j.$ Для определенности мы снова ограничимся случаем, когда параметр запаздывания $\tau$ меньше длины каждого ребра теперь уже дерева
${\cal T}.$ Однако аналогично построениям из раздела~7 в \cite{But23-RM}, можно рассмотреть и более общий случай -- например, когда
$2\tau<T,$ где $T$ -- высота дерева.

\newpage

\begin{center}
{\bf 4. Постановка вариационной задачи на произвольном дереве}
\end{center}

Рассмотрим компактное дерево ${\cal T}$ с множеством вершин $\{v_0,v_1,\ldots, v_m\}$ и множеством ребер $\{e_1,\ldots, e_m\}.$ Условимся,
что $\{v_0,v_{d+1},\ldots, v_m\}$ -- граничные вершины, т.е. принадлежащие только одному (граничному) ребру, а $\{v_1,\ldots, v_d\}$ --
внутренние.

Без ущерба для общности будем считать, что всякое ребро $e_j,\;j=\overline{1,m},$ начинается в вершине $v_{k_j}$ и заканчивается в $v_j,$ и
записывать $e_j=[v_{k_j},v_j],$ где $k_1=0.$ Вершину $v_0$ назовем {\it корнем}. Например, на рис.~2 имеем $m=9$ и $d=3,$ тогда как
$$
k_1=0, \quad k_2=k_3=1, \quad k_4=k_5=2, \quad k_6=k_7=k_8=k_9=3.
$$

Таким образом, $k_j$ порождает некоторое отображение $\{1,\ldots,m\}$ на $\{0,1,\ldots,d\},$ однозначно определяющее структуру ${\cal T}.$ А
именно, для всякого $j=\overline{0,d}$ множество $\{e_\nu\}_{\nu\in V_j},$ где $V_j:=\{\nu:k_\nu=j\},$ состоит из ребер, начинающихся в
вершине $v_j.$ В частности, имеем $\#\{e_\nu\}_{\nu\in V_0}=1,$ поскольку $v_0$ является граничной вершиной. Положим $k_j^{\{0\}}:=j$ и
$k_j^{\{\nu+1\}}:=k_{k_j^{\{\nu\}}}$ при $\nu=\overline{0,\nu_j},$ где $\nu_j$ таково, что $k_j^{\{\nu_j\}}=1.$ Тогда для каждого
$j=\overline{1,m}$ цепочка ребер $\{e_{k_j^{\{\nu\}}}\}_{\nu=\overline{0,\nu_j}}$ образует единственный простой путь между вершиной $v_j$ и
корнем. Обозначим через $T_j$ длину ребра $e_j.$ Величина $T:=\max_{j=\overline{d+1,m}}\sum_{\nu=0}^{\nu_j}T_{k_j^{\{\nu\}}}$ называется {\it
высотой} дерева~${\cal T}.$

\vspace*{-5mm}
\begin{center}
\unitlength=0.9mm
\begin{picture}(150,100)

 \put(75,18){\circle*{1}}         \put(73.5,13){\small $v_0$}
 \put(75,46){\circle*{1}}         \put(73.5,50){\small $v_1$}
 \put(51,61.95){\circle*{1}}      \put(44.5,64){\small $v_2$}
 \put(99,61.95){\circle*{1}}     \put(102,64){\small $v_3$}
 \put(25,49){\circle*{1}}         \put(19,46){\small $v_4$}
 \put(51,89.95){\circle*{1}}      \put(49.3,93.3){\small $v_5$}
 \put(79,82){\circle*{1}}         \put(74,84){\small $v_6$}
 \put(108,89){\circle*{1}}        \put(107,91.5){\small $v_7$}
 \put(127,61.95){\circle*{1}}     \put(129,61){\small $v_8$}
 \put(108,35){\circle*{1}}        \put(107.5,30.5){\small $v_9$}

 \put(75,18){\line(0,1){28}}      \put(70,31.5){\small $e_1$}
 \put(75,46){\line(-3,2){24}}     \put(60,50){\small $e_2$}
 \put(75,46){\line(3,2){24}}      \put(87,50){\small $e_3$}
 \put(51,61.95){\line(-2,-1){26}} \put(38,51.5){\small $e_4$}
 \put(51,61.95){\line(0,1){28}}   \put(52.5,76){\small $e_5$}
 \put(99,61.95){\line(-1,1){20}} \put(84.5,68.5){\small $e_6$}
 \put(99,61.95){\line(1,3){9}}   \put(99,78){\small $e_7$}
 \put(99,61.95){\line(1,0){28}}  \put(113.5,64){\small $e_8$}
 \put(99,61.95){\line(1,-3){9}}  \put(106,47){\small $e_9$}

 \put (30,0){\small Рис. 2. Пример нумерации вершин и ребер}

\end{picture}
\end{center}

\medskip

Пусть каждое ребро $e_j$ параметризовано переменной $t\in[0,T_j],$ причем $t=0$ соответствует его началу $v_{k_j},$ а $t=T_j$ -- концу $v_j.$
Под функцией $y$ на ${\cal T}$ будем понимать кортеж $y=[y_1,\ldots,y_m],$ чья $j\!$-ая компонента $y_j$ определена на ребре $e_j,$ т.е.
$y_j=y_j(t),\,t\in[0,T_j].$ Также зафиксируем $\tau\ge0$ и будем говорить, что функция $y$ определена на расширенном дереве ${\cal T}_\tau,$
если она определена на ${\cal T},$ а ее первая компонента $y_1(t)$ определена также при $t\in[-\tau,0).$

Пусть для определенности $\tau<T_j,\,j=\overline{1,m}.$ Рассмотрим управляемую систему, определяемую задачей Коши на ${\cal T}_\tau$ для
уравнений нейтрального типа:
\begin{equation}\label{2.1}
\ell_j y(t):=\sum_{k=0}^n\Big(b_{k,j}(t) y_j^{(k)}(t)+c_{k,j}(t)y_j^{(k)}(t-\tau)\Big)=u_j(t), \quad 0<t<T_j, \quad j=\overline{1,m},
\end{equation}
\begin{equation}\label{2.2}
y_j(t)=y_{k_j}(t+T_{k_j}), \quad t\in(-\tau,0), \quad j=\overline{2,m},
\end{equation}
\begin{equation}\label{2.2-1}
y_j^{(k)}(0)=y_{k_j}^{(k)}(T_j), \quad k=\overline{0,n-1}, \quad j=\overline{2,m},
\end{equation}
\begin{equation}\label{2.3}
y_1(t)=\varphi(t)\in W_2^n[-\tau,0], \quad t\in(-\tau,0),
\end{equation}
\begin{equation}\label{2.3-1}
y_1^{(k)}(0)=\varphi^{(k)}(0), \quad k=\overline{0,n-1},
\end{equation}
с комплекснозначными $\varphi(t)$ и $b_{k,j}(t),$ $c_{k,j}(t),$ $u_j(t),$ удовлетворяющими (\ref{1.8-1}).

Предполагается, что $j\!$-ое уравнение в (\ref{2.1}) определено на ребре $e_j$ дерева ${\cal T},$ причем (\ref{2.2-1}) являются условиями
склейки в его внутренних вершинах. Соотношения (\ref{2.2}) задают начальные функции для всех уравнений в (\ref{2.1}), кроме первого, и
означают, что запаздывание распространяется через все внутренние вершины дерева ${\cal T}.$ Предыстория процесса на всем дереве определяется
условиями (\ref{2.3}) и (\ref{2.3-1}).

Нетрудно показать, что задача Коши (\ref{2.1})--(\ref{2.3-1}) имеет единственное решение
$$
y=[y_1,\ldots,y_m]\in W_2^n({\cal T}_\tau):=W_2^n[-\tau,T_1]\oplus\bigoplus_{j=2}^m W_2^n[0,T_j],
$$
а при $d=1$ она совпадает с задачей (\ref{1.5})--(\ref{1.8}). Требуется найти управление
$$
u=[u_1,\ldots,u_m]\in L_2({\cal T}):=\bigoplus_{j=1}^m L_2(0,T_j),
$$
приводящее систему (\ref{2.1})--(\ref{2.3-1}) в состояние равновесия
\begin{equation}\label{2.4}
y_j(t)=0, \quad t\in[T_j-\tau,T_j], \quad j=\overline{d+1,m},
\end{equation}
и при этом минимизирующее норму $\|u\|_{L_2({\cal T})}=\sqrt{\sum_{j=1}^m\|u_j\|_{L_2(0,T_j)}^2}.$

Таким образом, приходим к вариационной задаче
\begin{equation}\label{2.5}
{\cal J}(y):=\sum_{j=1}^m\int_0^{T_j}|\ell_j y(t)|^2\,dt\to\min
\end{equation}
для функций $y=[y_1,\ldots,y_m]$ на ${\cal T}_\tau,$ удовлетворяющих условиям (\ref{2.2})--(\ref{2.4}).

Заметим, что условия (\ref{2.2}) никаких ограничений на функцию $y=[y_1,\ldots,y_m]$ не накладывают. Поэтому условимся, что взятие ${\cal
J}(y)$ и, в частности, $\ell_j y$ при $j=\overline{2,m}$ от какой бы то ни было функции $y$ на ${\cal T}$ автоматические подразумевает
применение условий (\ref{2.2}). Для краткости также введем обозначение $\ell y:=[\ell_1y,\ldots,\ell_my].$

\begin{center}
{\bf 5. Сведение к краевой задаче на дереве}
\end{center}

Рассмотрим в $W_2^k({\cal T}_\tau)$ обычное скалярное произведение
$$
(y,z)_{W_2^k({\cal T}_\tau)}=(y_1,z_1)_{W_2^k[-\tau,T_1]}+\sum_{j=2}^m (y_j,z_j)_{W_2^k[0,T_j]},
$$
где $y=[y_1,\ldots,y_m]$ и $z=[z_1,\ldots,z_m],$ тогда как $(f,g)_{W_2^k[a,b]}=\sum_{\nu=0}^k(f^{(\nu)},g^{(\nu)})_{L_2(a,b)}$ -- скалярное
произведение в $W_2^k[a,b],$ а $(\,\cdot\,,\,\cdot\,)_{L_2(a,b)}$ -- в $L_2(a,b).$

Обозначим через ${\cal W}$ замкнутое подпространство $W_2^n({\cal T}_\tau),$ состоящее из кортежей $[y_1,\ldots,y_m],$ удовлетворяющих
условиям (\ref{2.2-1}), (\ref{2.4}) и условию $y_1(t)=0$ на $[-\tau,0],$ причем если $\tau=0,$ то требуется $y_1^{(k)}(0)=y_j^{(k)}(T_j)=0,$
$j=\overline{d+1,m},$ $k=\overline{0,n-1}.$

Очевидно, на ${\cal W}$ можно смотреть и как на подпространство $W_2^n({\cal T})$ или $W_2^n(\widetilde{\cal T}),$ где $W_2^k({\cal
T}):=W_2^k({\cal T}_0),$ в то время как определение $W_2^k(\widetilde{\cal T})$ отличается от определения $W_2^k({\cal T}_0)$ только заменой
$T_j$ на $T_j-\tau$ при $j=\overline{d+1,m}.$

\medskip
{\bf Лемма 1. }{\it Пусть $y\in W_2^n({\cal T}_\tau)$ -- решение вариационной задачи (\ref{2.2})--(\ref{2.5}). Тогда
\begin{equation}\label{3.1}
B(y,w):=\sum_{j=1}^m\int_0^{T_j} \ell_j y(t)\overline{\ell_j w(t)}\,dt=0 \quad \forall \; w\in {\cal W}.
\end{equation}
Обратно, если $y\in W_2^n({\cal T}_\tau)$ удовлетворяет (\ref{2.2-1})--(\ref{2.4}) и (\ref{3.1}), то $y$ является решением задачи
(\ref{2.2})--(\ref{2.5}).}

\medskip
{\it Доказательство.} Пусть $y\in W_2^n({\cal T}_\tau)$ -- решение (\ref{2.2})--(\ref{2.5}). Тогда для всякого $w\in {\cal W}$ сумма $y+sw,$
в частности, при $s\in{\mathbb R}$ также принадлежит $W_2^n({\cal T}_\tau)$ и удовлетворяет условиям (\ref{2.2-1})--(\ref{2.4}). Положим
$F(s):={\cal J}(y+sw).$ Поскольку ${\cal J}(y+sw)\ge {\cal J}(y)$ для всех $s\in{\mathbb R},$ получаем $F'(0)=0.$ С другой стороны, имеет
место $F'(0)=2{\rm Re}B(y,w),$ что влечет (\ref{3.1}) в силу произвольности $w\in{\cal W}.$

Обратно, для всякой $y\in W_2^n({\cal T}_\tau),$ удовлетворяющей (\ref{2.2-1}), выполнение (\ref{3.1}) влечет
$$
{\cal J}(y+w)={\cal J}(y)+2{\rm Re}B(y,w)+{\cal J}(w)\ge {\cal J}(y)
$$
для всех $w\in {\cal W},$ что, в свою очередь, дает (\ref{2.5}) и при ограничениях (\ref{2.3})--(\ref{2.4}).  $\hfill\Box$

\medskip
Далее, применяя (\ref{2.2}) к $w=[w_1,\ldots,w_m]\in{\cal W},$ получаем представление
$$
\int_0^{T_j}\ell_jy(t)\overline{c_{k,j}(t)w_j^{(k)}(t-\tau)}\,dt= \int_0^{T_j-\tau}\ell_jy(t+\tau)\overline{c_{k,j}(t+\tau)w_j^{(k)}(t)}\,dt
\qquad\qquad\qquad\quad
$$
$$
\qquad\qquad\quad +\int_{T_{k_j}-\tau}^{T_{k_j}}\ell_jy(t+\tau-T_{k_j})\overline{c_{k,j}(t+\tau-T_{k_j})w_{k_j}^{(k)}(t)}\,dt, \;\;
j=\overline{1,m}, \;\; k=\overline{0,n},
$$
где $T_0=0$ и $w_0=0.$ Просуммировав по $j$ от $1$ до $m,$ будем иметь
$$
\sum_{j=1}^m \int_0^{T_j}\ell_jy(t)\overline{c_{k,j}(t)w_j^{(k)}(t-\tau)}\,dt= \sum_{j=1}^m
\int_0^{T_j-\tau}\ell_jy(t+\tau)\overline{c_{k,j}(t+\tau)w_j^{(k)}(t)}\,dt\qquad\qquad
$$
$$
\qquad\qquad\quad +\sum_{j=1}^d \sum_{\nu\in V_j} \int_{T_j-\tau}^{T_j}\ell_\nu
y(t+\tau-T_j)\overline{c_{k,\nu}(t+\tau-T_j)w_j^{(k)}(t)}\,dt, \quad k=\overline{0,n}.
$$

Таким образом, положив
$$
l_j:=T_j,\quad j=\overline{1,d}, \qquad l_j:=T_j-\tau, \quad j=\overline{d+1,m},
$$
в соответствии с (\ref{2.1}), перепишем (\ref{3.1}) в эквивалентном виде
\begin{equation}\label{3.2}
B(y,w)=\sum_{j=1}^m \sum_{k=0}^n \int_0^{l_j} \ell_{k,j} y(t) \overline{w_j^{(k)}(t)} \,dt=0\quad \forall \; w=[w_1,\ldots, w_n]\in {\cal W},
\end{equation}
где
\begin{equation}\label{3.3}
\ell_{k,j} y(t)=\overline{b_{k,j}(t)}\ell_j y(t) +\left\{\begin{array}{l}
\overline{c_{k,j}(t+\tau)}\ell_j y(t+\tau),\quad  0<t<T_j-\tau, \quad j=\overline{1,m},\\[4mm]
\displaystyle\sum_{\nu\in V_j}\overline{c_{k,\nu}(t+\tau-T_j)}\ell_\nu y(t+\tau-T_j),\\[-2mm]
          \qquad\qquad\qquad\qquad\qquad T_j-\tau< t<T_j, \quad j=\overline{1,d}.
\end{array}\right.
\end{equation}

В силу (\ref{1.8-1}) для всякой функции $y\in W_2^n({\cal T}_\tau)$ будем иметь
\begin{equation}\label{3.3-1}
\ell_{k,j} y(t)\in L(0,l_j), \quad k=\overline{0,n-1}, \quad \ell_{n,j} y(t)\in L_2(0,l_j), \quad j=\overline{1,m}.
\end{equation}

\medskip
{\bf Лемма 2. }{\it Пусть задан некоторый набор функций
$$
f_{k,j}(t)\in L(0,l_j), \quad k=\overline{0,n-1}, \quad f_{n,j}(t)\in L_2(0,l_j), \quad j=\overline{1,m}.
$$

1) Предположим, что имеет место тождество
\begin{equation}\label{3.4}
\sum_{j=1}^m \sum_{k=0}^n \int_0^{l_j} f_{k,j}(t) \overline{w_j^{(k)}(t)}\,dt=0  \quad \forall \; w=[w_1,\ldots, w_n]\in {\cal W}.
\end{equation}
Тогда
\begin{equation}\label{3.5-1}
g_{k,j}(t)\in W_1^1[0,l_j], \quad k=\overline{n,2n-1}, \quad j=\overline{1,m},
\end{equation}
и
\begin{equation}\label{3.6}
g_{k,j}(l_j)=\sum_{\nu\in V_j}g_{k,\nu}(0), \quad j=\overline{1,d}, \quad k=\overline{n,2n-1},
\end{equation}
а также
\begin{equation}\label{3.5-2}
g_{2n,j}(t)=0, \quad 0<t<l_j, \quad j=\overline{1,m},
\end{equation}
где $g_{k,j}(t)$ определяются формулами
\begin{equation}\label{3.6-1}
g_{n,j}(t):=f_{n,j}(t), \quad g_{n+l,j}(t):=f_{n-l,j}(t)-g_{n+l-1,j}'(t), \quad l=\overline{1,n}, \quad j=\overline{1,m}.
\end{equation}

2) Обратно, пусть выполняются (\ref{3.5-1}) и (\ref{3.6}). Тогда справедливо тождество
\begin{equation}\label{3.12-0}
\sum_{j=1}^m \sum_{k=0}^n \int_0^{l_j} f_{k,j}(t) \overline{w_j^{(k)}(t)}\,dt=\sum_{j=1}^m \int_0^{l_j} g_{2n,j}(t) \overline{w_j(t)}\,dt
\quad \forall \; w=[w_1,\ldots, w_n]\in {\cal W}.
\end{equation}}

\medskip

\medskip
{\it Доказательство.} 1) Выберем функции $F_{k,j}(t)\in W_1^{n-k}[0,l_j]$ так, чтобы
\begin{equation}\label{3.6-2}
F_{k,j}^{(n-k)}(t)=f_{k,j}(t), \quad k=\overline{0,n}, \quad j=\overline{1,m},
\end{equation}
и при этом удовлетворяющие условиям
\begin{equation}\label{3.7}
F_{k,j}^{(s)}(l_j)=\sum_{\nu\in V_j} F_{k,\nu}^{(s)}(0), \quad s=\overline{0,n-k-1}, \quad k=\overline{0,n-1}, \quad j=\overline{1,d},
\end{equation}
что, очевидно, возможно. Действительно, при $j=\overline{d+1,m}$ в качестве $F_{k,j}(t)$ можно выбрать любую первообразную $(n-k)\!$-го
порядка для $f_{k,j}(t),$ тогда как остальные первообразные $F_{k,j}(t)$ при $j=\overline{1,d}$ реккурентно определяются условиями
(\ref{3.7}).

Интегрируя по частям, получаем
\begin{equation}\label{3.8}
\begin{array}{c}
\displaystyle \sum_{j=1}^m \sum_{k=0}^n \int_0^{l_j} f_{k,j}(t) \overline{w_j^{(k)}(t)}\,dt
=\sum_{j=1}^m\sum_{k=0}^n (-1)^{n-k} \int_0^{l_j} F_{k,j}(t) \overline{w_j^{(n)}(t)}\,dt- \qquad\qquad\\[6mm]
\displaystyle \qquad\qquad\qquad\qquad\qquad\qquad -\sum_{j=1}^m\sum_{k=0}^{n-1}\sum_{s=1}^{n-k} (-1)^s
\Big(F_{k,j}^{(n-k-s)}(t)\overline{w_j^{(k+s-1)}(t)}\Big)\Big|_{t=0}^{l_j}.
\end{array}
\end{equation}
Кроме того, в силу определения класса ${\cal W},$ будем иметь
$$
\sum_{j=1}^m\Big(F_{k,j}^{(n-k-s)}(t)\overline{w_j^{(k+s-1)}(t)}\Big)\Big|_{t=0}^{l_j}=
$$
$$
= \sum_{j=1}^d F_{k,j}^{(n-k-s)}(l_j)\overline{w_j^{(k+s-1)}(l_j)}-\sum_{j=2}^m F_{k,j}^{(n-k-s)}(0)\overline{w_j^{(k+s-1)}(0)}=
$$
$$
= \sum_{j=1}^d \Big(F_{k,j}^{(n-k-s)}(l_j)-\sum_{\nu\in V_j} F_{k,\nu}^{(n-k-s)}(0)\Big)\overline{w_j^{(k+s-1)}(l_j)}=0, \;\;
k=\overline{0,n-1}, \;\; s=\overline{1,n-k},
$$
где последнее равенство следует из (\ref{3.7}). Таким образом, (\ref{3.4}) вместе с (\ref{3.8}) дает
\begin{equation}\label{3.9}
\sum_{j=1}^m\sum_{k=0}^n (-1)^{n-k} \int_0^{l_j} F_{k,j}(t) \overline{w_j^{(n)}(t)}\,dt =0, \quad w\in{\cal W},
\end{equation}
откуда следует, в частности, что
\begin{equation}\label{3.10}
\sum_{k=0}^n (-1)^{n-k} F_{k,j}(t)=\sum_{k=0}^{n-1}C_{k,j}\frac{t^k}{k!}, \quad j=\overline{1,m}.
\end{equation}
В самом деле, для этого достаточно ограничиться функциями $w_j(t)\in \!\!\stackrel{\rm o}{\,\,W_2^n}[0,\l_j]$ и учесть, что ортогональное
дополнение в $L_2(0,\l_j)$ множества их производных $n\!$-го порядка совпадает с множеством полиномов степени меньшей~$n.$

Поскольку $F_{k,j}(t)\in W_1^{n-k}[0,l_j]$ для $k=\overline{0,n-1}$ и $j=\overline{1,m},$ из (\ref{3.10}) вытекает
\begin{equation}\label{3.11}
{\cal F}_{l,j}(t):=\sum_{k=l}^n (-1)^{n-k} F_{k,j}(t)\in W_1^{n-l+1}[0,l_j], \quad l=\overline{0,n}, \quad j=\overline{1,m}.
\end{equation}
Используя (\ref{3.6-1}), (\ref{3.6-2}) и последнее обозначение, по индукции получаем
\begin{equation}\label{3.12}
g_{n+l,j}(t)=(-1)^l{\cal F}_{n-l,j}^{(l)}(t), \quad l=\overline{0,n}, \quad j=\overline{1,m},
\end{equation}
что вместе с (\ref{3.11}) дает (\ref{3.5-1}). Кроме того, из (\ref{3.10})--(\ref{3.12}) следует (\ref{3.5-2}).

Далее, подставляя (\ref{3.10}) в (\ref{3.9}) и интегрируя по частям, будем иметь
$$
0= \sum_{j=1}^m\sum_{k=0}^{n-1} C_{k,j} \int_0^{l_j}\frac{t^k}{k!} \overline{w_j^{(n)}(t)}\,dt= \qquad\qquad\qquad\qquad\qquad\qquad\quad
$$
$$
\,\;\;\;\qquad= \sum_{j=1}^m\sum_{k=0}^{n-1} C_{k,j} \Big(\sum_{s=0}^k(-1)^s\frac{l_j^{k-s}}{(k-s)!} \overline{w_j^{(n-s-1)}(l_j)} -(-1)^k
\overline{w_j^{(n-k-1)}(0)}\Big).
$$
Меняя порядок суммирования и учитывая определение ${\cal W},$ получаем
$$
\sum_{s=0}^{n-1}(-1)^s \sum_{j=1}^d \overline{w_j^{(n-s-1)}(l_j)} \Big(\sum_{k=s}^{n-1}C_{k,j}\frac{l_j^{k-s}}{(k-s)!} - \sum_{\nu\in V_j}
C_{s,\nu} \Big)=0.
$$
Отсутствие ограничений на величины $w_j^{(s)}(l_j)$ при $j=\overline{1,d}$ и $s=\overline{0,n-1}$ дает
\begin{equation}\label{3.11-1}
\sum_{k=s}^{n-1}C_{k,j}\frac{l_j^{k-s}}{(k-s)!} = \sum_{\nu\in V_j} C_{s,\nu}, \quad s=\overline{0,n-1}, \quad j=\overline{1,d}.
\end{equation}
Положим
\begin{equation}\label{3.11-0}
{\cal G}_{l,j}(t):=\sum_{k=0}^{l-1} (-1)^{n-k} F_{k,j}(t), \quad l=\overline{1,n}, \quad j=\overline{1,m},
\end{equation}
и запишем (\ref{3.10}) в виде
\begin{equation}\label{3.11-2}
{\cal G}_{l,j}(t)=\sum_{k=0}^{n-1}C_{k,j}\frac{t^k}{k!} -{\cal F}_{l,j}(t), \quad l=\overline{1,n}, \quad j=\overline{1,m},
\end{equation}
где ${\cal F}_{l,j}(t)$ определены в (\ref{3.11}). Согласно (\ref{3.7}) и (\ref{3.11-0}) будем иметь
\begin{equation}\label{3.11-3}
{\cal G}_{l,j}^{(s)}(l_j)=\sum_{\nu\in V_j} {\cal G}_{l,\nu}^{(s)}(0), \quad s=\overline{0,n-l}, \quad l=\overline{1,n}, \quad
j=\overline{1,d}.
\end{equation}
Подставляя (\ref{3.11-2}) в (\ref{3.11-3}) и используя (\ref{3.11-1}), приходим к равенствам
$$
\begin{array}{c}
\displaystyle {\cal F}_{l,j}^{(s)}(l_j) -\sum_{\nu\in V_j} {\cal F}_{l,\nu}^{(s)}(0) =\sum_{k=s}^{n-1}C_{k,j}\frac{l_j^{k-s}}{(k-s)!}
-\sum_{\nu\in V_j}
C_{s,\nu}=0, \qquad\qquad \\[5mm]
\qquad\qquad\qquad\qquad\qquad\qquad\qquad s=\overline{0,n-l}, \quad l=\overline{1,n}, \quad j=\overline{1,d}.
\end{array}
$$
Наконец, в силу (\ref{3.12}) последняя формула при $s=n-l$ дает (\ref{3.6}).

2) Предположим теперь, что имеют место (\ref{3.5-1}) и (\ref{3.6}). Интегрируя $n$ раз по частям с учетом (\ref{3.5-1}) и (\ref{3.6-1}),
получаем
$$
\sum_{k=0}^n \int_0^{l_j} f_{k,j}(t) \overline{w_j^{(k)}(t)} \,dt = \sum_{l=0}^{n-1}
\Big(g_{n+l,j}(t)\overline{w_j^{(n-l-1)}(t)}\Big)\Big|_{t=0}^{l_j} +\int_0^{l_j} g_{2n,j}(t)\overline{w_j(t)} \,dt.
$$
Кроме того, используя свойства функций класса ${\cal W}$ и (\ref{3.6}), будем иметь
$$
\sum_{j=1}^m\sum_{l=0}^{n-1} \Big(g_{n+l,j}(t)\overline{w_j^{(n-l-1)}(t)}\Big)\Big|_{t=0}^{l_j}=
\qquad\qquad\qquad\qquad\qquad\quad\qquad\qquad\qquad\qquad
$$
$$
=\sum_{l=0}^{n-1}\Big(\sum_{j=1}^d g_{n+l,j}(l_j)\overline{w_j^{(n-l-1)}(l_j)} -\sum_{j=2}^m g_{n+l,j}(0)\overline{w_j^{(n-l-1)}(0)}\Big)
$$
$$
\qquad\qquad\qquad\qquad\qquad\quad\;\; =\sum_{l=0}^{n-1}\sum_{j=1}^d \Big(g_{n+l,j}(l_j) -\sum_{\nu\in V_j}
g_{n+l,\nu}(0)\Big)\overline{w_j^{(n-l-1)}(l_j)}=0,
$$
что вместе с предыдущей формулой дает (\ref{3.12-0}). $\hfill\Box$

\medskip
Введем глобально нелокальные квазипроизводные $y_j^{\langle k\rangle}(t),$ $k=\overline{n,2n},$ $j=\overline{1,m},$ снова используя формулы
(\ref{1.8-2}) теперь вместе с (\ref{3.3}), и обозначим $y^{\langle k\rangle}:=[y_1^{\langle k\rangle},\ldots,y_m^{\langle k\rangle}].$

На функциях $y=[y_1,\ldots,y_m]\in W_2^n({\cal T}_\tau),$ удовлетворяющих условиям
\begin{equation}\label{3.11-6}
y_j^{\langle k\rangle}(t)\in W_1^1[0,l_j], \quad k=\overline{n,2n-1}, \quad j=\overline{1,m},
\end{equation}
рассмотрим краевую задачу, которую будем обозначать через ${\cal B},$ состоящую из функ\-ционально-дифференциальных уравнений $2n\!$-го
порядка
\begin{equation}\label{3.11-7}
y_j^{\langle 2n\rangle}(t)=0, \quad 0<t< l_j, \quad j=\overline{1,m},
\end{equation}
и условий (\ref{2.2})--(\ref{2.4}), а также нелокальных условий типа Кирхгофа
\begin{equation}\label{3.11-8}
y_j^{\langle k\rangle}(l_j)=\sum_{\nu\in V_j}y_\nu^{\langle k\rangle}(0), \quad j=\overline{1,d}, \quad k=\overline{n,2n-1}.
\end{equation}

Имеет место следующее утверждение.

\medskip
{\bf Лемма 3. }{\it Пусть $y\in W_2^n({\cal T}_\tau)$ удовлетворяет условиями (\ref{2.2-1}), (\ref{2.3}), (\ref{2.4}) и (\ref{3.1}). Тогда
имеет место (\ref{3.11-6}), а $y$ является решением краевой задачи  ${\cal B}.$

Обратно, всякое решение задачи ${\cal B}$ удовлетворяет (\ref{3.1}).}

\medskip
{\it Доказательство.} Принимая во внимание (\ref{3.3-1}), и применяя к (\ref{3.2}) первое утверждение леммы~2, приходим к
(\ref{3.11-6})--(\ref{3.11-8}).

Обратно, пусть $y\in W_2^n({\cal T}_\tau)$ является решением задачи ${\cal B}.$ Тогда второе утверждение леммы~2 вместе с левым равенством в
(\ref{3.2}), а также с (\ref{3.11-6}) и (\ref{3.11-8}) дает
\begin{equation}\label{3.12-2}
B(y,w)=\sum_{j=1}^m \int_0^{l_j} y_j^{\langle 2n\rangle}(t) \overline{w_j(t)}\,dt \quad \forall w\in {\cal W}.
\end{equation}
Наконец, (\ref{3.11-7}) влечет (\ref{3.1}). $\hfill\Box$

\medskip
{\bf Замечание 2.} В силу (\ref{3.12-2}) краевая задача ${\cal B}$ является самосопряженной в том смысле, что для всех $y,z\in{\cal W},$
удовлетворяющих условиям (\ref{3.11-6}) и (\ref{3.11-8}), имеет место соотношение $(y^{\langle 2n\rangle},z)_{L_2({\cal T})}=(y,z^{\langle
2n\rangle})_{L_2({\cal T})},$ где $(\,\cdot\,,\,\cdot\,)_{L_2({\cal T})}$ -- скалярное произведение в $L_2({\cal T}).$ Однако мы не будем
непосредственно ссылаться на этот факт.

\medskip
Объединяя леммы~1 и~3, приходим к основному результату настоящего раздела.

\medskip
{\bf Теорема 5. }{\it Функция $y\in W_2^n({\cal T}_\tau)$ является решением вариационной задачи (\ref{2.2})--(\ref{2.5}) тогда и только
тогда, когда она удовлетворяет (\ref{3.11-6}) и является решением краевой задачи ${\cal B}.$}

\begin{center}
{\bf 6. Однозначная разрешимость краевой задачи}
\end{center}

В данном разделе устанавливается однозначная разрешимость краевой задачи ${\cal B},$ а согласно теореме~5 -- и вариационной задачи
(\ref{2.2})--(\ref{2.5}). Введем обозначения
\begin{equation}\label{4.0-1}
\ell_j^0 y(t):=b_{n,j}(t)y_j^{(n)}(t)+c_{n,j}(t)y_j^{(n)}(t-\tau), \;\; \ell_j^1 y(t):=\ell_j y(t)-\ell_j^0 y(t), \;\;  j=\overline{1,m},
\end{equation}
$$
{\cal J}_\nu(y):= \sum_{j=1}^m\int_0^{T_j}|\ell_j^\nu y(t)|^2\,dt, \quad \nu=0,1,
$$
где автоматически предполагается (\ref{2.2}), а также $y_1(t)=0$ на $(-\tau,0).$

Рассмотрим банахово пространство $C^k(\widetilde{\cal T}),$ состоящее из функций $y=[y_1,\ldots,y_m],$ таких что $y_j(t)\in C^k[0,l_j],$
$j=\overline{1,m},$ с нормой $\|y\|_{C^k(\widetilde{\cal T})}:=\max_{j=\overline{1,m}}\|y_j\|_{C^k[0,l_j]}.$

При помощи теоремы Арцела нетрудно показать, что для каждого $k\in{\mathbb N}$ пространство $W_2^k(\widetilde{\cal T})$ компактно вложено в
$C^{k-1}(\widetilde{\cal T}),$ то есть всякое ограниченное множество в первом из них является предкомпактным во втором.

\medskip
{\bf Лемма 4. }{\it Существуют $C_0$ и $C_1,$ такие что
\begin{equation}\label{4.0}
{\cal J}(w)\le C_0\|w\|_{W_2^n(\widetilde{\cal T})}^2, \quad {\cal J}_1(w)\le C_1\|w\|_{C^{n-1}(\widetilde{\cal T})}^2 \quad \forall w\in
{\cal W}.
\end{equation}}

\medskip
{\it Доказательство.} Положим $b_j:=\|b_{n,j}\|_{L_\infty(0,T_j)}$ и $c_j:=\|c_{n,j}\|_{L_\infty(0,T_j)},$ а также
$$
\gamma:=\max_{j=\overline{1,m}\atop k=\overline{0,n-1}}\max\Big\{\|b_{k,j}\|_{L_2(0,T_j)},\|c_{k,j}\|_{L_2(0,T_j)}\Big\}.
$$
Пусть $w\in {\cal W}.$ Используя (\ref{2.1}) и неравенство
\begin{equation}\label{4.3-1}
(\alpha_1+\ldots+\alpha_s)^2\le s(\alpha_1^2+\ldots+\alpha_s^2), \quad \alpha_1,\ldots,\alpha_s\in{\mathbb R},
\end{equation}
для $s=2$ и $s=2n,$ получаем
$$
{\cal J}_0(w)\le 2\sum_{j=1}^m \Big(b_j^2\|w_j^{(n)}\|_{L_2(0,l_j)}^2+ c_j^2\int_0^{T_j}|w_j^{(n)}(t-\tau)|^2\,dt\Big),
$$
$$
{\cal J}_1(w)\le 2n\gamma^2\sum_{j=1}^m\sum_{k=0}^{n-1}\Big(\|w_j^{(k)}\|_{C[0,l_j]}^2+ \max_{0\le t\le T_j}|w_j^{(k)}(t-\tau)|^2\Big)
$$
соответственно. Учитывая (\ref{2.2}) применительно к $w,$ для всех $j=\overline{1,m}$ вычисляем
$$
\int_0^{T_j}|w_j^{(n)}(t-\tau)|^2\,dt=\|w_j^{(n)}\|_{L_2(0,T_j-\tau)}^2 + \|w_{k_j}^{(n)}\|_{L_2(T_{k_j}-\tau,T_{k_j})}^2,
$$
$$
\max_{0\le t\le T_j}|w_j^{(k)}(t-\tau)|= \max\big\{\|w_{k_j}^{(k)}\|_{C[T_{k_j}-\tau,T_{k_j}]}, \|w_j^{(k)}\|_{C[0,T_j-\tau]}\big\}, \quad
k=\overline{0,n-1},
$$
где $T_0=0$ и $w_0=0.$ Используя оценку ${\cal J}(w)\le2{\cal J}_0(w)+2{\cal J}_1(w),$ получаем (\ref{4.0}). $\hfill\Box$

\medskip
{\bf Лемма 5. }{\it Существует $c_0>0,$ такое что
\begin{equation}\label{4.1}
{\cal J}_0(w)\ge c_0\|w\|_{W_2^n(\widetilde{\cal T})}^2 \quad \forall w\in {\cal W}.
\end{equation}}

\medskip
{\it ДОказательство.} Согласно (\ref{4.0-1}), почти всюду на $(0,T_j)$ справедливы оценки
\begin{equation}\label{4.1-0}
|w_j^{(n)}(t)| \le \tilde b_j|\ell_j^0w(t)| +\tilde b_jc_j|w_j^{(n)}(t-\tau)|, \quad j=\overline{1,m},
\end{equation}
где $\tilde b_j:=\|b_{n,j}^{-1}\|_{L_\infty(0,T_j)},$ а $c_j$ определены в доказательстве предыдущей леммы.

Предположим от противного, что найдутся $w_{(s)}=[w_{(s),1},\ldots,w_{(s),m}]\in {\cal W},$ $s\in{\mathbb N},$ такие что
$\|w_{(s)}\|_{W_2^n(\widetilde{\cal T})}=1$ и
\begin{equation}\label{4.1-1}
{\cal J}_0(w_{(s)})\le \frac1s, \quad s\in{\mathbb N}.
\end{equation}
Используя (\ref{4.3-1}), (\ref{4.1-0}) и (\ref{4.1-1}), при всех $s\in{\mathbb N}$ и $j=\overline{1,m}$ получаем оценки
$$
\|w_{(s),j}^{(n)}\|_{L_2(l\tau,M_{l,j})} \le \frac{2\tilde b_j^2}s +2\tilde b_j^2c_j^2\|w_{(s),j}^{(n)}\|_{L_2((l-1)\tau,l\tau)},  \quad 0\le
l<\frac{l_j}\tau,
$$
где $M_{l,j}:=\min\{(l+1)\tau,l_j\}.$ В частности, будем иметь
$$
\lim_{s\to\infty}\|w_{(s),j}^{(n)}\|_{L_2(l\tau,M_{l,j})} =2\tilde b_j^2c_j^2\lim_{s\to\infty}\|w_{(s),j}^{(n)}\|_{L_2((l-1)\tau,l\tau)},
\quad 0\le l<\frac{l_j}\tau.
$$
Применяя тождества
$$
\|w_{(s),1}^{(n)}\|_{L_2(-\tau,0)}=0, \quad \|w_{(s),j}^{(n)}\|_{L_2(-\tau,0)}=\|w_{(s),k_j}^{(n)}\|_{L_2(T_{k_j}-\tau,T_{k_j})}, \quad
j=\overline{2,m}, \quad s\in{\mathbb N},
$$
по индукции получаем $\|w_{(s)}^{(n)}\|_{L_2({\cal T})}\to0$ при $s\to\infty,$ где $w^{(n)}=[w_1^{(n)},\ldots, w_m^{(n)}].$

С другой стороны, аналогично лемме~5 в \cite{But23-arXiv} нетрудно показать, что $\|w^{(n)}\|_{L_2({\cal T})}$ порождает норму в ${\cal W},$
эквивалентную норме $\|w\|_{W_2^n(\widetilde{\cal T})},$ а значит, и $\|w_{(s)}\|_{W_2^n(\widetilde{\cal T})}\to0$ при $s\to\infty.$
Последнее противоречит сделанному предположению. $\hfill\Box$

\medskip
{\bf Лемма 6. }{\it Существует $c>0,$ такое что
$$
{\cal J}(w)\ge c\|w\|_{W_2^n(\widetilde{\cal T})}^2 \quad \forall w\in {\cal W}.
$$}

\medskip
{\it Доказательство.} Снова предположим от противного, что найдутся $w_{(s)}\in {\cal W}$ при $s\in{\mathbb N},$ такие что
$\|w_{(s)}\|_{W_2^n(\widetilde{\cal T})}=1,$ но теперь
\begin{equation}\label{4.3}
{\cal J}(w_{(s)})\le \frac1s, \quad s\in{\mathbb N}.
\end{equation}
Неравенство ${\cal J}_0(w)\le2{\cal J}(w)+2{\cal J}_1(w)$ вместе с оценками (\ref{4.0}) и (\ref{4.1}) влечет
\begin{equation}\label{4.5}
\frac{c_0}2\|w\|_{W_2^n(\widetilde{\cal T})}^2\le {\cal J}(w) +C_1\|w\|_{C^{n-1}(\widetilde{\cal T})}^2, \quad w\in{\cal W}.
\end{equation}

В силу компактности вложения $W_2^n(\widetilde{\cal T})$ в $C^{n-1}(\widetilde{\cal T})$ найдется подпоследовательность $\{w_{(s_k)}\},$
сходящаяся в $C^{n-1}(\widetilde{\cal T}).$ Неравенство (\ref{4.5}) дает
$$
\frac{c_0}2\|w_{(s_k)}-w_{(s_l)}\|_{W_2^n(\widetilde{\cal T})}^2 \le {\cal J}(w_{(s_k)}-w_{(s_l)})+
C_1\|w_{(s_k)}-w_{(s_l)}\|_{C^{n-1}(\widetilde{\cal T})}^2.
$$
Кроме того, в силу (\ref{4.3}) имеем ${\cal J}(w_{(s_k)}-w_{(s_l)})\le2/s_k+2/s_l,$ а значит, подпоследовательность $\{w_{(s_k)}\}$ является
фундаментальной в ${\cal W}$ и имеет там предел $w_{(0)}.$

В силу  леммы~4 сходимость $w_{(s_k)}$ к $w_{(0)}$ в ${\cal W}$ влечет $\ell w_{(s_k)} \to \ell w_{(0)}$ в $L_2({\cal T}).$ Следовательно, в
силу (\ref{4.3}), будем иметь
$$
\|\ell w_{(0)}\|_{L_2({\cal T})}^2= \lim_{k\to\infty}\|\ell w_{(s_k)}\|_{L_2({\cal T})}^2=\lim_{k\to\infty}{\cal J}(w_{(s_k)})=0,
$$
то есть $\ell w_{(0)}=0.$ Итак, $w_{(0)}$ является решением задачи Коши (\ref{2.1})--(\ref{2.3-1}) с $\varphi(t)\equiv 0$ и $u_j(t)=0$ при
$j=\overline{1,m},$ а значит, $w_{(0)}=0,$ что противоречит $\|w_{(0)}\|_{W_2^n(\widetilde{\cal T})}=1.$ $\hfill\Box$

\medskip
Следующая теорема является основным результатом данного раздела.

\medskip
{\bf Теорема 6. }{\it Краевая задача ${\cal B}$ имеет единственное решение $y\in W_2^n({\cal T}_\tau),$ удовлетворяющее условиям
(\ref{3.11-6}). Кроме того, выполняется оценка
\begin{equation}\label{4.8}
\|y\|_{W_2^n({\cal T}_\tau)}\le C\|\varphi\|_{W_2^n[-\tau,0]},
\end{equation}
где $C$ не зависит от $\varphi(t).$}

\medskip
{\it Доказательство.} Рассмотрим функцию $\Phi=[\Phi_1,\ldots,\Phi_m]\in W_2^n({\cal T}_\tau)$ вида
$$
\Phi_1(t)=\left\{\begin{array}{cc}
\varphi(t), & -\tau\le t<0,\\[1mm]
\varphi_1(t),& 0\le t\le T_1-\tau,\\[1mm]
0, & T_1-\tau<t\le T_1,
\end{array}\right.
\qquad \Phi_j(t)\equiv0, \quad j=\overline{2,m},
$$
где
$$
\varphi_1(t)=\sum_{k=0}^{n-1}\varphi^{(k)}(0)p_k(t),
$$
а $p_k(t)$ являются базисными многочленами Эрмита степени $2n-1,$ удовлетворяющие условиям $p_k^{(\nu)}(0)=\delta_{k,\nu}$ и
$p_k^{(\nu)}(T_1-\tau)=0$ при $\nu,k=\overline{0,n-1}.$ Здесь $\delta_{k,\nu}$ -- символ Кронекера.

В силу леммы~3 всякая функция $y\in W_2^n({\cal T}_\tau),$ для которой выполняются условия (\ref{2.2-1}), (\ref{2.3}) и (\ref{2.4}), будет
удовлетворять условиями абсолютной непрерывности квазипроизводных (\ref{3.11-6}) и являться решением краевой задачи ${\cal B}$ тогда и только
тогда, когда имеет место (\ref{3.1}).

Другими словами, $y$ удовлетворяет (\ref{3.11-6}) и является решением задачи ${\cal B},$ если и только если $x:=y-\Phi\in {\cal W}$ и
\begin{equation}\label{4.9}
B(\Phi,w)+B(x,w)=0 \quad \forall w\in {\cal W}.
\end{equation}

Так как $B(w,w)={\cal J}(w),$ полуторалинейная форма $(\,\cdot\,,\,\cdot\,)_{\cal W}:=B(\,\cdot\,,\,\cdot\,),$  в силу леммы~6, является
скалярным произведением в ${\cal W}.$ Кроме того, справедлива оценка
$$
|B(\Phi,w)|=\Big|\int_0^{T_1}\ell_1\Phi(t)\overline{\ell_1w(t)}\,dt\Big| \le \qquad\qquad\qquad\qquad\qquad\quad\;
$$
\begin{equation}\label{4.10}
\qquad\qquad\qquad\qquad\qquad\le M\|\varphi\|_{W_2^n[-\tau,0]}\|w\|_{W_2^n(\widetilde{\cal T})}\le
\frac{M}{\sqrt{c}}\|\varphi\|_{W_2^n[-\tau,0]}\|w\|_{\cal W},
\end{equation}
где $\|w\|_{\cal W}:=\sqrt{(w,w\,)_{\cal W}}.$ Таким образом, в силу теоремы Рисса об общем виде линейного ограниченного функционала в
гильбертовом пространстве, существует единственная функция $x\in {\cal W},$ такая что выполняется равенство (\ref{4.9}), а значит, задача
${\cal B}$ имеет единственное решение $y=\Phi+x.$ Согласно (\ref{4.9}) и (\ref{4.10}) имеем
$$
\|x\|_{\cal W}\le \frac{M}{\sqrt{c}}\|\varphi\|_{W_2^n[-\tau,0]}.
$$
Снова применяя лемму~6, приходим к (\ref{4.8}). $\hfill\Box$

\begin{center}
{\bf 7. Сравнение введенных квазипроизводных с классическими}
\end{center}

Различные формы квазипроизводных играют важную роль в спектральной теории обыкновенных дифференциальных операторов (см., например,
\cite{Neum, Shin, Krein, Glaz, NZadeSchk, SavchShk, Vlad-04,12-HryPro, MirShk-16, Vlad17-arXiv, Bond-22, Bond-23}). Принципиальным отличием
введенных квазипроизводных (\ref{1.4-0}) является нелокальность. Однако в локальном случае $\tau=0$ их можно сравнить с известными. Для этой
цели приведем некоторые сведения о квазипроизводных. Рассмотрим дифференциальное выражение четного порядка
\begin{equation}\label{7.1}
L_ny:=\sum_{k,s=0}^n (r_{ks}(t)y^{(n-k)})^{(n-s)}, \quad -\infty\le a<t<b\le\infty,
\end{equation}
понимаемое обычным образом, коль скоро каждый коэффициент $r_{ks}(t)$ имеет производные до порядка $n-s$ включительно. Чтобы придать смысл
выражению (\ref{7.1}), например, при локально суммируемых $r_{ks}(t),$ можно ввести квазипроизводные
$$
y^{\{n\}}:=\sum_{k=0}^n r_{k0}(t)y^{(n-k)}, \quad y^{\{n+l\}}:=(y^{\{n+l-1\}})'+\sum_{k=0}^n r_{kl}(t)y^{(n-k)}, \quad l=\overline{1,n},
$$
предполагая $y^{\{k\}}\in AC_{\rm loc}(a,b)$ при $k=\overline{n,2n-1}.$ При этом, очевидно, $L_ny=y^{\{2n\}}.$

Тогда выражение (\ref{7.1}) можно записать в эквивалентном векторном виде
\begin{equation}\label{7.2}
{\cal L}Y=Q_0(t)Y'-Q(t)Y
\end{equation}
с локально суммируемыми матрицами $Q_0(t)$ и $Q(t)$ размерности $2n\times2n,$ где
$$
Y=[y,y',\ldots,y^{(n-1)},y^{\{n\}},\ldots,y^{\{2n-1\}}]^T, \quad {\cal L}Y=[0,\ldots,0,L_ny]^T,
$$
причем свободный член для (\ref{7.2}) должен иметь нули на тех же местах, что и ${\cal L}Y.$

Чтобы матрица $Q_0(t)$ была почти всюду обратимой, а элементы $(Q_0(t))^{-1}$ принадлежали $L_{\infty,\rm loc}(a,b),$ достаточно потребовать
$r_{00}^{-1},r_{0k}\in L_{\infty,\rm loc}(a,b)$ при $k=\overline{1,n}.$

Таким образом квазипроизводные использовались, например, в \cite{Krein, Glaz, Neum}, но для выражений более частного вида, нежели
(\ref{7.1}). Вообще говоря, квазипроизводные можно задавать по-разному, но все такие построения укладываются в общую конструкцию
квазидифференциального выражения, предложенную Д.~Шином \cite{Shin}.

Обобщением указанной выше цели применения квазипроизводных является регуляризация сингулярных дифференциальных выражений с коэффициентами из
пространств обобщенных функций \cite{NZadeSchk, SavchShk, Vlad-04, 12-HryPro, MirShk-16, Vlad17-arXiv, Bond-22, Bond-23}, то есть
преобразование таких дифферен\-циальных выражений к эквивалентному виду в пространстве вектор-функций, но с регулярными (суммируемыми либо
локально суммируемыми) коэффициентами.

В соответствии с этим различные формы квазипроизводных можно ранжировать по степени сингулярности коэффициентов дифференциального выражения,
которое они регуляризуют \cite{Bond-22}. Для выражений четного порядка наиболее ''сильные'' в таком смысле квазипроизводные введены
К.А.~Мирзоевым и А.А.~Шкаликовым \cite{MirShk-16}.

Регуляризация в \cite{MirShk-16} допускает для коэффициентов в (\ref{7.1}) выполнение условий
\begin{equation}\label{sing}
|r_{00}|^{-1/2},|r_{00}|^{-1/2}R_{ks}\in L_{2,\rm loc}(a,b), \quad k+s=\overline{1,2n},
\end{equation}
где $R_{ks}$ -- первообразная порядка $\min\{k,s\}$ для $r_{ks}$ в смысле распределений.

В частности, это соответствует регуляризации выражения Штурма--Лиувилля
\begin{equation}\label{7.2-1}
Ly:=-y''+ q(t)y, \quad 0<t<1,
\end{equation}
с потенциалом-распределением $q\in W_2^{-1}[0,1]$ (см. \cite{SavchShk, NZadeSchk}). Последнее означает, что $q=\sigma'$ для некоторой функции
$\sigma\in L_2(0,1).$ Тогда квазипроизводная
\begin{equation}\label{7.3}
y^{[1]}=y'-\sigma y
\end{equation}
приводит к соответствующему квазидифференциальному выражению
$$
y^{[2]}:=(y^{[1]})'+\sigma y^{[1]}+\sigma^2y=-Ly.
$$

Отметим, что (в случае локальной суммируемости) матрица $(Q_0(t))^{-1}Q(t)$ в (\ref{7.2}) называется {\it ассоциированной} с дифференциальным
выражением (\ref{7.1}) и принадлежит в терминологии \cite{MirShk-16} классу Шина--Зеттла. Своя ассоциированная матрица из этого класса
возникает при любом способе регуляризации и хранит информацию о форме использованных квазипроизводных. Например, для (\ref{7.3}) она имеет
вид
\begin{equation}\label{7.4-1}
\left[\begin{array}{cc}
\sigma & 1\\[0mm]
-\sigma^2 & -\sigma
\end{array}\right].
\end{equation}
Обратно, всякая матрица из класса Шина--Зеттла (см. определение в \cite{MirShk-16}) порождает некоторое квазидифференциальное выражение
соответствующего порядка.

Теперь посмотрим, как вписываются в данную картину квазипроизводные (\ref{1.4-0}). Пусть $n=1.$ При $\tau=0,\,T=1,$ $b_1=1$ и $c_0=c_1=0$
выражение в (\ref{1.1}) примет вид
$$
\ell y=y'+by, \quad b:=b_0\in L_2(0,1).
$$
Тогда, согласно (\ref{1.4-00}), имеем $\tilde\ell_0y=\overline{b}y'+|b|^2y,$ $\tilde\ell_1y=y'+by,$ что вместе с (\ref{1.4-0}) дает
\begin{equation}\label{q-10}
y^{\langle 1\rangle}=y'+by,
\end{equation}
$$
y^{\langle 2\rangle}=\tilde\ell_0 y -(y^{\langle 1\rangle})' =-y'' -i(({\rm Im}\,b)y)' -i({\rm Im}\,b)y' +(|b|^2-({\rm Re}\,b)')y.
$$
В частности, для вещественных $b$ будем иметь $y^{\langle 2\rangle}=Ly$ вида (\ref{7.2-1}), где
\begin{equation}\label{q-11-2}
q=-b'+b^2,
\end{equation}
что может быть достигнуто и без требования вещественности $b,$ но путем отказа от комплексного сопряжения в (\ref{1.4-00}). Тогда
(\ref{7.2-1}) представляется в виде
$$
y^{\langle 2\rangle}=-(y^{\langle 1\rangle})' +by^{\langle 1\rangle}.
$$

Таким образом, квазипроизводная (\ref{q-10}) регуляризует выражение Штурма--Лиу\-вил\-ля (\ref{7.2-1}) с потенциалами $q\in W_2^{-1}[0,1],$
допускающими представление (\ref{q-11-2}). Вместе с тем вещественность $q,$ как и при использовании квазипроизводной (\ref{7.3}), очевидно,
не требуется. Соответствующая ассоциированная матрица имеет вид
\begin{equation}\label{q-11-4}
\left[\begin{array}{cc}
-b & 1\\[0mm]
0 & b
\end{array}\right].
\end{equation}
При этом класс допустимых $q\in W_2^{-1}[0,1]$ определяется разрешимостью в $L_2(0,1)$ уравнения Риккати (\ref{q-11-2}), которое приводится к
уравнению Штурма--Лиувилля
\begin{equation}\label{q-14}
-u''+qu=0
\end{equation}
заменой $b=-u'/u.$ Чтобы $b\in L_2(0,1),$ решение $u$ уравнения (\ref{q-14}) не должно обращаться в ноль на отрезке $[0,1].$ Потенциалы $q,$
для которых такое решение существует и вещественно, называются потенциалами Миуры \cite{12-HryPro}. Для всякого комплексного $q\in
W_2^{-1}[0,1]$ существование такого решения гарантировано после прибавления дос\-таточно большого положительного числа $C,$ то есть заменой
$q$ на $q+C.$ В частности, вещественный потенциал $q\in W_2^{-1}[0,1]$ является потенциалом Миуры, если оператор Дирихле, порожденный левой
частью (\ref{q-14}), положительно определен \cite{12-HryPro}. Согласно лемме~6 из предыдущего раздела, это условие является и необходимым.

Итак, для потенциалов $q\in W_2^{-1}[0,1]$ квазипроизводная $y^{\langle1\rangle}$ совпадает по ''силе'' в указанном выше смысле с
квазипроизводной $y^{[1]}$ и (с учетом модификации для комплексных~$q)$ дает альтернативную регуляризацию выражения (\ref{7.2-1}) после
сдвига спектрального параметра при необходимости. Кроме того, $y^{\langle1\rangle}$ совпадает с квазипроизводной в \cite{12-HryPro}, где
требовалась факторизация
$$
-\frac{d^2}{dt^2}+q=-\Big(\frac{d}{dt}-b\Big)\Big(\frac{d}{dt}+b\Big).
$$
В \cite{Beals-88} использовалась аналогичная факторизация выражений произвольного порядка на оси, но не с целью регуляризации или введения
квазипроизводных.

Преемственность матриц (\ref{7.4-1}) и (\ref{q-11-4}) прослеживается в одном результате из \cite{Bond-23}, согласно которому все
ассоциированные с выражением (\ref{7.2-1}) матрицы имеют вид
$$
\left[\begin{array}{cc}
\sigma_1 & 1\\[0mm]
q_1-\sigma_1^2 & -\sigma_1
\end{array}\right],
\quad \sigma_1\in L_2(0,1), \quad  q_1\in L(0,1), \quad q=\sigma_1' +q_1.
$$
В \cite{Bond-23} также получено описание ассоциированных матриц для нормального вида (см. \cite{MirShk-16}) общего выражения (\ref{7.1}) в
случае $r_{00}=1.$

Пусть теперь $n$ произвольное. Рассмотрим соответствующее локальное выражение $\ell$ в (\ref{1.1}) вместе с формально самосопряженным ему
выражением $\ell^*:$

$$
\ell y=\sum_{k=0}^n b_k(t)y^{(k)}, \quad \ell^* y=\sum_{k=0}^n (-1)^k (\overline{b_k(t)}y)^{(k)}, \quad 0<t<1,
$$
где коэффициенты удовлетворяют условиям (\ref{1.1-1}) при $a=1.$ Тогда (\ref{1.4-00}), (\ref{1.4-0}) дают, в частности,
$y^{\langle2n\rangle}=\ell^*\ell y,$ что проще увидеть, сопоставляя формулы (\ref{3.1}) и (\ref{3.12-2}). При этом соответствующие
квазипроизводные (\ref{1.4-0}) регуляризуют выражение $\ell^*\ell.$

Согласно лемме~6 оператор $\ell^*\ell y$ с краевыми условиями
\begin{equation}\label{q-16}
y^{(k)}(0)=y^{(k)}(1)=0, \quad k=\overline{0,n-1},
\end{equation}
положительно определен. Рассмотрим дифференциальное выражение $L_n$ вида (\ref{7.1}) на интервале $(0,1)$ и вместо (\ref{sing}) для простоты
предположим
$$
r_{00},r_{00}^{-1}\in L_\infty(0,1), \quad r_{ks} \in W_2^{-l}[0,1], \quad l:=\min\{k,s\}, \quad k+s=\overline{1,2n}.
$$
В свете сказанного выше представляет интерес вопрос: Допускает ли это выражение факторизацию $L_n=\ell^*\ell,$ если оператор $L_n y$ с
краевыми условиями (\ref{q-16}) положительно определен либо полуограничен снизу с достаточно большой константой?

В случае положительного ответа можно было бы снова говорить об альтернативной \cite{MirShk-16} регуляризации сингулярного выражения
(\ref{7.1}), но когда оно является формально самосопряженным. Как было продемонстрировано выше, для самосопряженного сингулярного выражения
Штурма--Лиувилля $L$ эта гипотеза верна.

\medskip
Автор выражает благодарность А.Л.~Скубачевскому за рекомендацию посмотреть задачу об успокоении системы управления с последействием в связи с
предложенной автором идеей глобального запаздывания на графе, а также Н.П.~Бондаренко, М.Ю.~Игнатьеву, М.А.~Кузнецовой и П.А.~Терехину за
полезные обсуждения.

\begin{center}
{\bf Литература}
\end{center}

\end{document}